\documentclass[10pt,reqno]{amsart}

\usepackage{amssymb}
\usepackage{xypic}
\xyoption{all}
\def\umono{\ar@{_{(}->}[u]}
\def\uumono{\ar@{_{(}->}[uu]}

\def\lmono{\ar@{_{(}->}[l]}
\def\llmono{\ar@{_{(}->}[ll]}

\newcommand{\Z}{{\mathbb Z}}

\newcommand{\C}{{\mathbb C}}
\newcommand{\F}{{\mathbb F}}

\newcommand{\tor}{\operatorname{Tor}\nolimits}

\newcommand{\A}{\ifmmode{\mathcal{A}}\else${\mathcal{A}}$\fi}
\newcommand{\K}{\ifmmode{\mathcal{K}}\else${\mathcal{K}}$\fi}
\newcommand{\U}{\ifmmode{\mathcal{U}}\else${\mathcal{U}}$\fi}
\newcommand{\T}{\ifmmode{\mathcal{T}}\else${\mathcal{T}}$\fi}
\newcommand{\FF}{\ifmmode{\mathcal{F}}\else${\mathcal{F}}$\fi}
\newcommand{\LL}{\ifmmode{\mathcal{L}}\else${\mathcal{L}}$\fi}

\newtheorem{theorem}{Theorem}[section]
\newtheorem{proposition}[theorem]{Proposition}
\newtheorem{corollary}[theorem]{Corollary}
\newtheorem{lemma}[theorem]{Lemma}

\title[Homotopy exponents for large $H$-spaces] {Homotopy exponents for large $H$-spaces}

\author{Wojciech Chach\'olski}
\address{K.T.H.\\ Matematik\\ S-10044 Stockholm,
Sweden} \email{wojtek@math.kth.se}

\author{Wolfgang Pitsch}
\author{J\'er\^ome Scherer}
\address{Universitat Aut\`onoma de Barcelona \\
Departament de Mate\-m\`ati\-ques\\
E-08193 Bellaterra, Spain} \email{pitsch@mat.uab.es,
jscherer@mat.uab.es}

\author{Don Stanley}
\address{Department of Mathematics and Statistics\\
University of Regina, College West 307.14\\
Regina, Saskatchewan\\
Canada  S4S 0A2} \email{Donald.Stanley@uregina.ca}

\thanks{The first author is supported in
part by NSF grant DMS-0296117, Vetenskapsr\aa det grant 2001-4296,
and G\"oran Gustafssons Stiftelse. The second and third authors
are partially supported by FEDER/MEC grant MTM2007-61545. The
second author is supported by the program Ram\'on y Cajal, MEC,
Spain.}


\begin{document}


\begin{abstract}
We show that $H$-spaces with finitely generated cohomology, as an
algebra or as an algebra over the Steenrod algebra, have homotopy
exponents at all primes. This provides a positive answer to a
question of Stanley.
\end{abstract}


\maketitle


\section*{Introduction}
\label{sec intro}

A simply connected space is elliptic if both its rational homotopy
and rational homology are finite. Moore's conjecture, see for
example \cite{MR952582}, predicts that elliptic complexes have an
exponent at any prime $p$, meaning that there is a bound on the
$p$-torsion in the graded group of all homotopy groups. Any finite
$H$-space is known to be elliptic as it is rationally equivalent
to a finite product of (odd dimensional) spheres. Relying on
results by James \cite{MR0079263} and Toda \cite{MR0092968} about
the homotopy groups of spheres, the fourth author (re)proved in
\cite{MR1888808} Long's result that finite $H$-spaces have an
exponent at any prime \cite{Long}. He proved in fact a stronger
result which holds for example for $H$-spaces for which the mod
$p$ cohomology is finite. He also asked whether this would hold
for finitely generated cohomology rings. The aim of this note is
to give a positive answer to this question and provide a way
larger class of $H$-spaces which have homotopy exponents.

\medskip

\noindent {\bf Theorem \ref{thm main}} {\it Let $X$ be a connected
and $p$-complete $H$-space such that $H^*(X; \F_p)$ is finitely
generated as an algebra over the Steenrod algebra. Then $X$ has an
exponent at $p$.}

\medskip

This class of $H$-spaces is optimal in the sense that $H$-spaces
with a larger mod $p$ cohomology, such as an infinite product of
Eilenberg-Mac Lane spaces $K(\Z/p^n, n)$, will not have in general
an exponent at~$p$. As a corollary, we obtain the desired result.
In fact we obtain the following global theorem.

\medskip

\noindent {\bf Theorem \ref{thm fg}} {\it Let $X$ be a connected
$H$-space such that $H^*(X; \Z)$ is finitely generated as an
algebra. Then $X$ has an exponent at each prime $p$.}

\medskip

The methods we use are based on the deconstruction techniques of
the third author in his joint work with Castellana and Crespo,
\cite{MR2264802}. Our results on homotopy exponents should also be
compared with the computations of homological exponents done with
Cl\'ement, \cite{H*exponents}. Whereas such $H$-spaces always have
homotopy exponents, they almost never have homological exponents.
The only simply connected $H$-spaces for which the $2$-torsion in
$H_*(X; \Z)$ has a bound are products of mod $2$ finite $H$-spaces
with copies of the infinite complex projective space $\C P^\infty$
and $K(\Z, 3)$.

\medskip

{\bf Acknowledgments.} This project originated during a workshop
at the CRM, Barcelona, held in the emphasis year on algebraic
topology (2007-08). We would like to thank the organizers for
making it possible to meet in such a pleasant atmosphere.

\section{Homotopy exponents}
\label{sec Serre}

Our starting point is the fact that mod $p$ finite $H$-spaces have
always homotopy exponents. The following is a variant of Stanley's
\cite[Corollary~2.9]{MR1888808}. Whereas he focused on spaces
localized at a prime, we will stick to $p$-completion in the sense
of Bousfield and Kan, \cite{MR51:1825}. Since the $p$-localization
map $X \rightarrow X_{(p)}$ is a mod $p$ homology equivalence, his
result implies the following.

\begin{proposition}[Stanley]
\label{Don}
Let $p$ be a prime and $X$ be a $p$-complete and connected
$H$-space such that $H^*(X; \F_p)$ is finite. Then $X$ has an
exponent at $p$.
\end{proposition}

We will not repeat the proof, but let us sketch the main steps.
Let us consider a decomposition of $X$ by $p$-complete cells, i.e.
$X$ is obtained by attaching cones along maps from
$(S^n)^\wedge_p$. The natural map $X \rightarrow \Omega \Sigma X$
factors then through the loop spaces on a wedge $W$ of a finite
numbers of such $p$-completed spheres, up to multiplying by some
integer $N$: the composite $X \rightarrow \Omega \Sigma X
\xrightarrow{N} \Omega \Sigma X$ is homotopic to $X \rightarrow
\Omega W \rightarrow \Omega \Sigma X$. The proof goes by induction
on the number of $p$-complete cells and the key ingredient here is
Hilton's description of the loop space on a wedge of spheres,
\cite{MR0068218}. Note that the suspension of a map between
spheres is torsion except for the multiples of the identity. This
idea to ``split off" all the cells of $X$ up to multiplication by
some integer is dual to Arlettaz' way to split off Eilenberg-Mac
Lane spaces in $H$-spaces with finite order $k$-invariants,
\cite[Section~7]{MR1775748}. The final step relies on the
classical results by James, \cite{MR0079263}, and Toda,
\cite{MR0092968}, that spheres do have homotopy exponents at all
primes.

\begin{theorem}
\label{thm main}
Let $X$ be a connected and $p$-complete $H$-space such that
$H^*(X; \F_p)$ is finitely generated as an algebra over the
Steenrod algebra. Then $X$ has an exponent at $p$.
\end{theorem}

\begin{proof}
A connected $H$-space such that $H^*(X; \F_p)$ is finitely
generated as an algebra over the Steenrod algebra can always be
seen as the total space of an $H$-fibration $F \rightarrow X
\rightarrow Y$ where $Y$ is an $H$-space with finite mod $p$
cohomology and $F$ is a $p$-torsion Postnikov piece whose homotopy
groups are finite direct sums of copies of cyclic groups $\Z/p^r$
and Pr\"ufer groups $\Z_{p^\infty}$,
\cite[Theorem~7.3]{MR2264802}. This is a fibration of $H$-spaces
and $H$-maps, so that we obtain another fibration $F^\wedge_p
\rightarrow X \rightarrow Y^\wedge_p$ by $p$-completing it. The
base space $Y^\wedge_p$ now satisfies the assumptions of
Proposition~\ref{Don}. It has therefore an exponent at $p$. The
homotopy groups of the fiber $F^\wedge_p$ are finite direct sums
of cyclic groups $\Z/p^n$ and copies of the $p$-adic integers
$\Z^\wedge_p$. Thus $F^\wedge_p$ has an exponent at $p$ as well.
The homotopy long exact sequence of the fibration allows us to
conclude.
\end{proof}

We see here how the $p$-completeness assumption plays an important
role. The space $K(\Z_{p^\infty}, 1)$ for example has obviously no
exponent at $p$, but its $p$-completion is $K(\Z^\wedge_p, 2) =
(\C P^\infty)^\wedge_p$, which is a torsion free space. The mod
$p$ cohomology of $K(\Z_{p^\infty}, 1)$ is a polynomial ring on
one generator in degree~$2$, we must thus also work with
$p$-complete spaces to give an answer to Stanley's question
\cite[Question~2.10]{MR1888808}.

\begin{corollary}
\label{cor main}
Let $X$ be a connected and $p$-complete $H$-space such that
$H^*(X; \F_p)$ is finitely generated as an algebra. Then $X$ has
an exponent at $p$. \hfill{\qed}
\end{corollary}

In fact, when the mod $p$ cohomology is finitely generated, the
fiber $F$ in the fibration described in the proof of
Theorem~\ref{thm main} is a single Eilenberg-Mac Lane space $K(P,
1)$. Thus the typical example of an $H$-space with finitely
generated mod $p$ cohomology is the $3$-connected cover of a
simply connected finite $H$-space ($P$ is $\Z_{p^\infty}$ in this
case). Likewise, the typical example in Theorem~\ref{thm main} are
highly connected covers of finite $H$-spaces. This explains why
such spaces have homotopy exponents!

If one does not wish to work at one prime at a time and prefers to
find a global condition which permits to conclude that a certain
class of spaces have exponents at all primes, one must replace mod
$p$ cohomology by integral cohomology.

\begin{theorem}
\label{thm fg}
Let $X$ be a connected $H$-space such that $H^*(X; \Z)$ is
finitely generated as an algebra. Then $X$ has an exponent at each prime $p$.
\end{theorem}

\begin{proof}
Since the integral cohomology groups are finitely generated it
follows from the universal coefficient exact sequence (see
\cite{MR592883}) that the integral homology groups are also
finitely generated. Since $X$ is an $H$-space we may use a
standard Serre class argument to conclude that so are the homotopy
groups. Therefore  the $p$-completion map $X \rightarrow
X^\wedge_p$ induces an isomorphim on the $p$-torsion at the level
of homotopy groups. The theorem is now a direct consequence of the
next lemma.
\end{proof}

\begin{lemma}
\label{lemma comparefg}
Let $X$ be a connected space. If $H^*(X; \Z)$ is finitely
generated as an algebra, then so is $H^*(X; \F_p)$.
\end{lemma}

\begin{proof}
Let $u_1, \dots , u_r$ generate $H^*(X; \Z)$ as an algebra.
Consider the universal coefficients short exact sequences
$$
0 \rightarrow H^n(X; \Z) \otimes \Z/p \longrightarrow H^n(X; \F_p)
\xrightarrow{\partial} \tor (H^{n+1}(X; \Z); \Z/p) \rightarrow 0\,
.
$$
Since $H^{*}(X; \Z)$ is finitely generated as an algebra it is
degree-wise finitely generated as a group and therefore $\tor
(H^{*}(X; \Z); \Z/p)$ can be identified with the ideal of elements
of order $p$ in  $H^*(X; \Z)$. This ideal must be finitely
generated since $H^*(X; \Z)$ is Noetherian. Choose generators
$a_1, \dots, a_s$. Each $a_i$ corresponds to a pair $\alpha_i,
\beta \alpha_i$ in $H^*(X; \F_p)$, where $\beta$ denotes the
Bockstein.

We claim that the elements $\alpha_1, \dots , \alpha_s$ together
with the mod $p$ reduction of the algebra generators, denoted by
$\bar u_1, \dots, \bar u_r$, generate $H^*(X; \F_p)$ as an
algebra. Let $x \in H^*(X; \F_p)$ and write its image $\partial(x)
= \sum \lambda_j a_j$ with $\lambda_j= \lambda_j(u)$ a polynomial
in the $u_i$'s. Define now $\bar \lambda_j = \lambda_j(\bar u) \in
H^*(X; \F_p)$ to be the corresponding polynomial in the $\bar
u_i$'s. As the action of $H^*(X; \Z)$ on the ideal $\tor (H^{*}(X;
\Z); \Z/p)$ factors through the mod $p$  reduction map $H^*(X; \Z)
\rightarrow H^*(X; \F_p)$, the element $x - \sum \bar \lambda_j
\alpha_j$ belongs to the kernel of $\partial$, i.e. it lives in
the image of the mod $p$ reduction. It can be written therefore as
a polynomial $\bar \mu$ in the $\bar u_i$'s. Thus $x = \bar \mu +
\sum \bar\lambda_j \alpha_j$.
\end{proof}

\providecommand{\bysame}{\leavevmode\hbox
to3em{\hrulefill}\thinspace}
\providecommand{\MR}{\relax\ifhmode\unskip\space\fi MR }
\providecommand{\MRhref}[2]{%
  \href{http://www.ams.org/mathscinet-getitem?mr=#1}{#2}
} \providecommand{\href}[2]{#2}

\end{document}